\documentclass[twoside,a4paper,11pt]{amsart}

\usepackage[margin=0.95in]{geometry}

\usepackage{amsfonts, amsmath,amssymb}
\usepackage[utf8]{inputenc}
\usepackage{hyperref}
\usepackage[all]{xy}
\usepackage[lite]{amsrefs}
\usepackage{graphicx}
\usepackage{wrapfig}
\usepackage{enumerate}
\usepackage[british]{babel}
\usepackage{tikz}
\usetikzlibrary{matrix,arrows}
\newcommand*{\Steinberg}{\widetilde{\operatorname{H}}}

\newcommand{\hyper}{\mathcal{H}}
\newcommand{\calO}{\mathcal{O}}
 \newcommand{\SLtwo}{\mathrm{SL}_2}
 
 \newcommand{\SLO}{\mathrm{SL}_2\left(\mathcal{O}_{-m}\right)}
 \newcommand{\PSLO}{\mathrm{PSL}_2\left(\mathcal{O}_{-m}\right)}
 
\newcommand{\Q}{\mathbb Q}

\newcommand{\Dsix}{{\mathcal{D}}_6}
\newcommand{\Dfour}{{\mathcal{D}}_4}
\newcommand{\Dthree}{{\mathcal{D}}_3}
\newcommand{\Dtwo}{{\mathcal{D}}_2}
\newcommand{\Dn}{\mathcal{D}_n}
\newcommand{\Dtwon}{\mathcal{D}_{n}}
\newcommand{\Dl}{\mathcal{D}_\ell}
\newcommand{\Af}{{\mathcal{A}}_4}
\newcommand{\Afour}{{\mathcal{A}}_4}

\newcommand{\Sfour}{{\mathcal{S}}_4}
\newcommand{\Farrell}{\widehat{\operatorname{H}}}
\newcommand{\C}{{\mathbb{C}}}
\newcommand{\N}{{\mathbb{N}}}
\newcommand{\symmetricspace}{\mathcal{H}}
\newcommand{\ra}{\ensuremath{\rightarrow}}
\newcommand{\PSL}{\mathrm{PSL}}
\newcommand{\sign}{{\rm sign}}
\newcommand*{\Cohomol}{\operatorname{H}}
\newcommand{\ef}{{\mathbb F}_2}
\newcommand{\xsp}{X_s^\prime}
\newcommand{\F}{{\mathbb{F}}}
\newcommand{\Z}{{\mathbb{Z}}}
\newcommand{\R}{{\mathbb{R}}}
\newcommand{\CC}{\ensuremath{\mathbb{C}}}
\newcommand{\rationals}{{\mathbb{Q}}}
\newcommand{\ringO}{\mathcal{O}}
\newcommand{\ringOm}{\mathcal{O}_{-m}}
\newcommand{\Hy}{\mathcal{H}}

\newcommand*{\Homol}{\operatorname{H}}
\newcommand{\op}{\operatorname}

\newcommand{\imQuadRing}{{\mathcal{O}}_{\Q(\sqrt{-m}\thinspace )}}
\newcommand{\arithGrp}{{\rm SL}_2({\mathbb{Z}}[\sqrt{-2}\thinspace][\frac{1}{2}])}

\theoremstyle{plain}
\newtheorem{thm}{\bfseries Theorem}
\newtheorem{lemma}[thm]{\bfseries Lemma}
\newtheorem{theorem}[thm]{\bfseries Theorem}
\newtheorem{proposition}[thm]{\bfseries Proposition}
\newtheorem{conjecture}[thm]{\bfseries Conjecture}
\newtheorem{corollary}[thm]{\bfseries Corollary}
\theoremstyle{definition}
\newtheorem{df}[thm]{\bfseries Definition}
\newtheorem{remark}[thm]{\bfseries Remark}
\newtheorem*{isomorphismCondition}{\bfseries Condition B}
\newtheorem*{ConditionBprime}{\bfseries Condition B'}

\newtheorem*{ConditionA}{\bfseries Condition A}

\newcommand{\cellCondition}{A}
\newcommand{\weakerCondition}{B}
\newcommand{\IsomorphismConditionSymbol}{B}
\newcommand{\technicalCondition}{B'}

\usepackage{pstricks} 

\newcommand{\circlegraph}{ 
\begin{pspicture}     (-0.21,-0.155)(0.21,0.155) 
             \pscircle(0,0.0){0.15}
             \psdots(0.15,0.0)
\end{pspicture} }
 
\newcommand{\edgegraph}{ 
\begin{pspicture}(-0.3,-0.1)(0.3,0.3)
\psdots(-0.2,0.0)
\psdots(0.2,0.0)
\psline(-0.2,0.0)(0.2,0.0)
\end{pspicture} }

\newcommand{\graphFive}{  
\begin{pspicture}(-0.21,-0.155)(0.21,0.155)
\pscircle(0,0.0){0.15}
\psdots(-0.15,0)
\psdots(0.15,0)
\psline(-0.15,0)(0.15,0)
\end{pspicture} }

\newcommand{\graphTwo}{  
\begin{pspicture}(0,0.05)(0.8,0.4)
\pscircle(0.2,0.2){0.15}
\psdots(0.35,0.2)
\psline(0.35,0.2)(0.65,0.2)
\psdots(0.65,0.2)
\end{pspicture} }

\begin{document}
\pagestyle{plain}

\title{An Introduction to Torsion Subcomplex Reduction}
  \date{\today}
  \author{Alexander D. Rahm}
  \address{Laboratoire de math\'ematiques GAATI, Universit\'e de la Polyn\'esie Fran\c{c}aise, BP 6570 -- 98702 Faaa, French Polynesia}
\urladdr{http://gaati.org/rahm/}
\email{Alexander.Rahm@upf.pf}

\maketitle

 \begin{abstract}
  This survey paper introduces to a technique called Torsion Subcomplex Reduction (TSR) for computing torsion in the cohomology of discrete groups acting on suitable cell complexes.
TSR enables one to skip machine computations on cell complexes, and to access directly the reduced torsion subcomplexes,
which yields results on the cohomology of matrix groups in terms of formulas.
  TSR has already yielded general formulas for the cohomology of the tetrahedral Coxeter groups as well as, at odd torsion, of SL$_2$ groups over arbitrary number rings. The latter formulas allow to refine the Quillen conjecture.
Furthermore, progress has been made to adapt TSR to Bredon homology computations.
In particular for the Bianchi groups, yielding their equivariant $K$-homology, and,
by the Baum--Connes assembly map, the $K$-theory of their reduced $C^*$-algebras.
As a side application, TSR has allowed to provide dimension formulas for the Chen--Ruan orbifold cohomology of the complexified Bianchi orbifolds, and to prove Ruan's crepant resolution conjecture for all complexified Bianchi orbifolds.
 \end{abstract}

\section{Introduction}
This survey paper is based on the habilitation thesis of the author, restricting to the expository parts, which are updated here,
and referring to previously published papers for the proofs.
The goal is to introduce to a technique for computing Farrell--Tate cohomology of arithmetic groups, presented in Section~\ref{techniques}.
This technique can also be applied in the computation of other invariants, as described in Section~\ref{Results}, where further results are stated.

\subsection{Background}\label{background}

Our objects of study are discrete groups~$\Gamma$ such that~$\Gamma$ admits a torsion-free subgroup of finite index.
By a theorem of Serre~\cite{SerreGroupesDiscrets}, all the torsion-free subgroups of finite index in~$\Gamma$ have the same cohomological dimension;
this dimension is called the virtual cohomological dimension (abbreviated vcd) of~$\Gamma$.
Above the vcd, the (co)homology of a discrete group is determined by its system of finite subgroups.
We are going to discuss it in terms of Farrell--Tate cohomology.
The Farrell--Tate cohomology $\Farrell^q$  is  identical to group cohomology $\Homol^q$ in all degrees $q$ above the vcd,
and extends in lower degrees to a cohomology theory of the system of finite subgroups.
Details are elaborated in Brown's book~\cite{Brown}*{chapter X}.
So for instance considering the Coxeter groups, the virtual cohomological dimension of all of which vanishes,
their Farrell--Tate cohomology is identical to all of their group cohomology.
In Section~\ref{conjugacy reduction}, we will introduce a method of how to explicitly determine the Farrell--Tate cohomology:
By reducing torsion sub-complexes.

Let us note that for the same arithmetic groups, cohomology outside of our setting has much stronger contemporary interest,
and therefore, there has been extensive work on it. Just to mention a few, fairly recent publications about group cohomology in low cohomological degrees, from which to find more references:
On SL$_N(\Z)$ with rising rank $N$ and modulo small torsion~\cites{sikiri2019voronoi},
on infinite towers of congruence subgroups \cites{AGMY,BergeronSengunVenkatesh},
on arbitrary groups using general purpose algorithms \cites{Ellis}.

\subsection{Overview of the results}

This paper introduces the technique of \emph{torsion subcomplex reduction}.
It is a technique for the study of discrete groups $\Gamma$,
giving easier access to the cohomology of the latter at a fixed prime~$\ell$
and above the virtual cohomological dimension,
by extracting the relevant portion of the equivariant spectral sequence and then simplifying it.
Instead of having to work with a full cellular complex $X$ with a nice $\Gamma$-action,
the technique inputs only an often lower-dimensional subcomplex of $X$,
and reduces it to a small number of cells.

The author first developed torsion subcomplex reduction for a specific class of arithmetic groups,
the Bianchi groups, for which the method yielded all of the homology above the virtual cohomological dimension~\cite{Rahm:homological_torsion}.
Some elements of this technique had already been used by Soul\'e for a modular group~\cite{Soule};
and were used by Mislin and Henn as a set of ad hoc tricks.
After rediscovering these ad hoc tricks, the author puts them into a general framework~\cite{Rahm:formulas}.
The advantage of using this framework
 is that it becomes possible to find general formulas for the dimensions of the Farrell--Tate cohomology, for instance for the entire family of the Bianchi groups.

It is convenient to give some examples of where the technique of torsion subcomplex reduction has already produced good results:
\begin{itemize}
 \item The Bianchi groups and their congruence subgroups (cf. Section \ref{The Bianchi groups});
 \item The Coxeter groups (cf. Section \ref{The Coxeter groups});
 \item The SL$_2$ groups over arbitrary number rings (cf. Section \ref{formulas for the Farrell--Tate cohomology});
 \item PSL$_4(\Z)$ and the PGL\texorpdfstring{$_3$}{(3)} groups over rings of quadratic integers (cf. Section \ref{GL3}).
 \item The technique has also been adapted to groups with non-trivial centre (cf. Section \ref{non-trivial-centre}).
\end{itemize}
This has led to the following applications:
\begin{itemize}
 \item  Refining the Quillen conjecture (cf. Section \ref{QC}),
 \item Computing equivariant \textit{K}-homology (cf. Section \ref{Bredon state}),
 \item Understanding Chen--Ruan orbifold cohomology (cf. Section \ref{orbifold state}).
\end{itemize}

\section{The technique of Torsion Subcomplex Reduction} \label{techniques}

\subsection{Farrell--Tate cohomology and Steinberg homology}\label{Farrell--Tate cohomology and Steinberg homology}
Let $\Gamma$ be a virtual duality group: this means,
$\Gamma$ admits a finite index subgroup $\Gamma'$ such that $\Z$ admits a finite projective resolution over $\Z[\Gamma']$,
and there is an integer $n$ such that $\Homol^i(\Gamma; \thinspace \Z[\Gamma])= 0$ for $i\neq n$ and $\Homol^n(\Gamma; \thinspace \Z[\Gamma])$ is $\Z$-torsion-free.
Then $\Gamma$ is of finite virtual cohomological dimension vcd$(\Gamma) = n < \infty$ with $n$ the aforementioned integer
(where we have to make the smallest choice $n=0$ if $n$ is not unique).
Then the ``dualizing module'' is $D := \Homol^n(\Gamma;\thinspace \Z[\Gamma])$,
and the \textit{Steinberg homology} of $\Gamma$
(with coefficients $M$)
is $\Steinberg_i(\Gamma; \thinspace M) := \Homol_i(\Gamma; \thinspace D\otimes M)$.
Recall~\cite{Brown79}*{\S 11.8} that there is an exact sequence
tying together group cohomology $\Homol^\bullet$, Steinberg homology $\Steinberg_\bullet$ and Farrell-Tate cohomology $\Farrell^\bullet$ of $\Gamma$ :
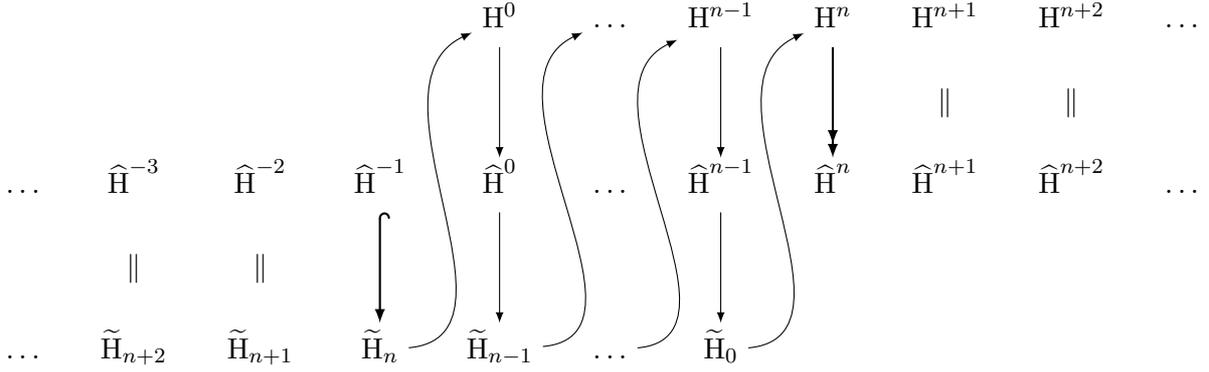
\begin{center}
\begin{tikzpicture}[descr/.style={fill=white,inner sep=1.5pt}]
        \matrix (m) [
            matrix of math nodes,
            row sep=1em,
            column sep=1.4em,
            text height=1.99ex, text depth=0.75ex
        ]
        {       &                 &                &              & \Homol^0       & \hdots & \Homol^{n-1} & \Homol^{n} & \Homol^{n+1} & \Homol^{n+2} & \hdots \\
	        &                 &                &              &                &        &              &            & \parallel    & \parallel    &        \\
         \hdots & \Farrell^{-3}   & \Farrell^{-2}  & \Farrell^{-1}& \Farrell^{0}   & \hdots &\Farrell^{n-1}&\Farrell^{n}&\Farrell^{n+1}&\Farrell^{n+2}& \hdots \\
                & \parallel       & \parallel      &              &                &        &              &            &              &              &        \\
         \hdots &\Steinberg_{n+2}&\Steinberg_{n+1}&\Steinberg_{n}&\Steinberg_{n-1}& \hdots &\Steinberg_{0}\\
        };

        \path[overlay,->, font=\scriptsize,>=latex]
         (m-5-4) edge[out=-355,in=-155] (m-1-5)
         (m-3-4) edge[thick, right hook->] (m-5-4)
         (m-1-5) edge (m-3-5)
         (m-3-5) edge (m-5-5)
         (m-1-7) edge (m-3-7)
         (m-3-7) edge (m-5-7)
         (m-5-5) edge[out=-355,in=-155] (m-1-6)
         (m-5-6) edge[out=-355,in=-155] (m-1-7)
         (m-5-7) edge[out=-355,in=-155] (m-1-8)
         (m-1-8) edge[thick, draw,->>] (m-3-8)
;
\end{tikzpicture}
\end{center}
Therefore, Brown describes the Farrell--Tate cohomology of $\Gamma$ to consist of the cohomology functors $\Homol^i$ for $i>n$,
the Steinberg homology functors $\Steinberg_i$ for $i>n$,
modified $\Homol^n$ and $\Steinberg_n$ functors, and $n$ additional functors $\Farrell^0$, $\hdots$, $\Farrell^{n-1}$
which are some sort of mixture of the functors $\Homol^i$ and $\Steinberg_i$ for $i \leq n$.

\subsection{Reduction of torsion subcomplexes in the classical setting} \label{conjugacy reduction}
Let $\ell$ be a prime number.
We require any discrete group $\Gamma$
under our study to be provided with what we will call a \textit{polytopal $\Gamma$-cell complex}, that is,
a finite-dimensional simplicial complex $X$ with cellular
$\Gamma$-action such that each cell stabiliser fixes its cell point-wise.
In practice, we relax the simplicial condition to a polytopal one,
merging finitely many simplices to a suitable polytope.
We could obtain the simplicial complex back as a triangulation.
We further require that the fixed point set~$X^G$ be acyclic for every non-trivial finite $\ell$-subgroup $G$ of~$\Gamma$.

Then, the $\Gamma$-equivariant Farrell--Tate cohomology $\Farrell^*_\Gamma(X; \thinspace M)$ of~$X$, for any trivial $\Gamma$-module $M$ of coefficients, gives us the
$\ell$-primary part $\Farrell^*(\Gamma; \thinspace M)_{(\ell)}$ of the Farrell--Tate cohomology of~$\Gamma$, as follows.
\begin{proposition}[Brown \cite{Brown}] \label{Brown's proposition}
For a $\Gamma$-action on $X$ as specified above, the canonical map
$$ \Farrell^*(\Gamma; \thinspace M)_{(\ell)} \to \Farrell^*_\Gamma(X; \thinspace M)_{(\ell)} $$
is an isomorphism.
\end{proposition}

The classical choice \cite{Brown} is to take for $X$
 the geometric realization of the partially ordered set of non-trivial finite subgroups
(respectively, non-trivial elementary Abelian $\ell$-subgroups) of~$\Gamma$,
 the latter acting by conjugation. The stabilisers are then the normalizers, which in many discrete groups are infinite.
In addition, there are often great computational challenges to determine a group presentation for the normalizers.
When we want to compute the module $\Farrell^*_\Gamma(X; \thinspace M)_{(\ell)}$
subject to Proposition~\ref{Brown's proposition},
at least we must know the ($\ell$-primary part of the) Farrell--Tate cohomology of these normalizers.
The Bianchi groups are an instance where different isomorphism types can occur for this cohomology
 at different conjugacy classes of elementary Abelian $\ell$-subgroups, both for $\ell=2$ and $\ell=3$.
As the only non-trivial elementary Abelian $3$-subgroups in the Bianchi groups are of rank $1$,
the orbit space $_\Gamma \backslash X$ consists only of one point for each conjugacy class of type $\Z/3$
and a corollary~\cite{Brown} from Proposition~\ref{Brown's proposition} decomposes the
 $3$-primary part of the Farrell--Tate cohomology of the Bianchi groups into the direct product over their normalizers.
However, due to the different possible homological types of the normalizers (in fact, two of them occur),
 the final result remains unclear and subject to tedious case-by-case computations of the normalizers.

In contrast, in the cell complex we are going to construct (specified in Definition~\ref{reduced torsion subcomplex definition} below),
 the connected components of the orbit space are for the $3$-torsion in the Bianchi groups not simple points,
 but have either the shape $\edgegraph$ or $\circlegraph$.
This dichotomy already contains the information about the occurring normalizer.

The starting point for our construction is the following definition.

\begin{df}
 Let $\ell$ be a prime number. The \emph{$\ell$-torsion subcomplex} of a polytopal $\Gamma$-cell complex~$X$
 consists of all the cells of $X$ whose stabilisers in~$\Gamma$  contain elements of order $\ell$.
\end{df}

We are from now on going to require the cell complex $X$ to admit only finite stabilisers in~$\Gamma$,
and we require the action of $\Gamma$ on the coefficient module $M$ to be trivial.
Then obviously only cells from the \emph{$\ell$-torsion subcomplex} contribute to $\Farrell^*_\Gamma(X; \thinspace M)_{(\ell)}$.

\begin{corollary}[Corollary to Proposition~\ref{Brown's proposition}] \label{Brownian}
 There is an isomorphism between the $\ell$-primary parts of the Farrell--Tate cohomology of~$\Gamma$ and the
 $\Gamma$-equivariant Farrell--Tate cohomology of the $\ell$-torsion subcomplex.
\end{corollary}

We are going to reduce the \emph{$\ell$-torsion subcomplex} to one which still carries the
$\Gamma$-equivariant Farrell--Tate cohomology of~$X$,
but which can also have considerably fewer orbits of cells.
This can be easier to handle in practice,
and, for certain classes of groups, leads us to an explicit structural description of the Farrell--Tate cohomology of~$\Gamma$.
The pivotal property of this reduced $\ell$-torsion subcomplex will be given in Theorem~\ref{pivotal}.
Our reduction process uses the following conditions,
which are imposed to a triple $(\sigma, \tau_1, \tau_2)$
of cells in the $\ell$-torsion subcomplex,
where $\sigma$ is a cell of dimension $n-1$,
lying in the boundary of precisely the two $n$-cells $\tau_1$ and~$\tau_2$,
 the latter cells representing two different orbits.

\begin{ConditionA} \label{cell condition}
The triple $(\sigma, \tau_1, \tau_2)$ is said to satisfy Condition A
if no higher-dimensional cells of the $\ell$-torsion subcomplex touch $\sigma$,
if the interior of $\tau_1$ and the interior of $\tau_2$ do not contain two points which are on the same orbit,
and if the $n$-cell stabilisers admit an isomorphism
$\Gamma_{\tau_1} \cong \Gamma_{\tau_2}$.
\end{ConditionA}

Where this condition is fulfilled in the $\ell$-torsion subcomplex,
 we merge the cells $\tau_1$ and $\tau_2$ along~$\sigma$ and do so for their entire orbits,
 if and only if they meet the following additional condition B.
We will refer by \emph{mod $\ell$ cohomology}
to group cohomology with $\Z/\ell$-coefficients under the trivial action.

\begin{isomorphismCondition}
With the notation above Condition $\cellCondition$, the inclusion $ \Gamma_{\tau_1} \subset \Gamma_\sigma$ induces an isomorphism on mod $\ell$ cohomology.
\end{isomorphismCondition}

\begin{lemma}[\cite{Rahm:formulas}] \label{A}
 Let $\widetilde{X_{(\ell)}}$ be the $\Gamma$-complex obtained by orbit-wise merging two $n$-cells of the
 $\ell$-torsion subcomplex $X_{(\ell)}$
which satisfy Conditions~$\cellCondition$ and~$\IsomorphismConditionSymbol$.
Then, $$\Farrell^*_\Gamma(\widetilde{X_{(\ell)}}; \thinspace M)_{(\ell)} \cong \Farrell^*_\Gamma(X_{(\ell)}; \thinspace M)_{(\ell)}.$$
\end{lemma}

By a ``terminal $(n-1)$-cell'',
we will denote an $(n-1)$-cell $\sigma$ with
\begin{itemize}
 \item modulo~$\Gamma$
precisely one adjacent $n$-cell $\tau$,
\item and such that $\tau$ has no further cells on the $\Gamma$-orbit of $\sigma$ in its boundary;
\item there shall be no higher-dimensional cells adjacent to $\sigma$.
\end{itemize}
And by ``cutting off'' the $n$-cell $\tau$,
 we will mean that we remove $\tau$ together with~$\sigma$ from our cell complex.

\begin{df} \label{reduced torsion subcomplex definition}
 A \emph{reduced $\ell$-torsion subcomplex} associated to a polytopal $\Gamma$-cell complex~$X$
 is a cell complex obtained by recursively merging orbit-wise all the pairs of cells satisfying
 conditions~$\cellCondition$ and~$\weakerCondition$,
 and cutting off $n$-cells that admit a terminal $(n-1)$-cell when condition~$\weakerCondition$ is satisfied.
\end{df}

A priori, this process yields a unique reduced $\ell$-torsion subcomplex only up to suitable isomorphisms,
so we do not speak of ``the'' reduced $\ell$-torsion subcomplex.
The following theorem makes sure that the $\Gamma$-equivariant mod $\ell$ Farrell--Tate cohomology
is not affected by this issue.

\begin{theorem}[\cite{Rahm:formulas}] \label{pivotal}
 There is an isomorphism between the $\ell$-primary part of the Farrell--Tate cohomology of~$\Gamma$ and the
 $\Gamma$-equivariant Farrell--Tate cohomology of a reduced $\ell$-torsion subcomplex obtained from $X$ as specified above.
\end{theorem}

In order to have a practical criterion for checking Condition~$\IsomorphismConditionSymbol$,
we make use of the following stronger condition.

Here, we write ${\rm N}_{\Gamma_\sigma}$ for taking the normalizer in ${\Gamma_\sigma}$ and
${\rm Sylow}_\ell$ for picking an arbitrary Sylow $\ell$-subgroup.
 This is well defined because all Sylow $\ell$-subgroups are conjugate.
We use Zassenhaus's notion for a finite group to be $\ell$-\emph{normal},
 if the center of one of its Sylow $\ell$-subgroups is the center of every Sylow $\ell$-subgroup in which it is contained.

\begin{ConditionBprime}
With the notation of Condition $\cellCondition$, the group $\Gamma_\sigma$ admits a (possibly trivial) normal subgroup $T_\sigma$ with trivial mod~$\ell$ cohomology
and with quotient group $G_\sigma$; and the group $\Gamma_{\tau_1}$ admits a (possibly trivial) normal subgroup
$T_\tau$ with trivial mod~$\ell$ cohomology and with quotient group $G_\tau$ making the sequences
\begin{center}
 $ 1 \to T_\sigma \to \Gamma_\sigma \to G_\sigma \to 1$ and $ 1 \to T_\tau \to \Gamma_{\tau_1} \to G_\tau \to 1$
\end{center}
exact and satisfying one of the following.
\begin{enumerate}
 \item  Either $G_\tau \cong G_\sigma$, or
 \item $G_\sigma$ is $\ell$-normal and $G_\tau \cong {\rm N}_{G_\sigma}({\rm center}({\rm Sylow}_\ell(G_\sigma)))$, or
 \item both $G_\sigma$  and $G_\tau$ are $\ell$-normal and there is a (possibly trivial) group $T$
 with trivial mod~$\ell$ cohomology making the sequence
$$1 \to T \to {\rm N}_{G_\sigma}({\rm center}({\rm Sylow}_\ell(G_\sigma))) \to {\rm N}_{G_\tau}({\rm center}({\rm Sylow}_\ell(G_\tau))) \to 1$$
exact.
\end{enumerate}
\end{ConditionBprime}

\begin{lemma}[\cite{Rahm:formulas}] \label{Implying the isomorphism condition}
Condition B' implies Condition B.
\end{lemma}

\begin{remark}
 The computer implementation \cite{BuiRahm:scpInHAP} checks Conditions~$\technicalCondition (1)$ and $\technicalCondition (2)$ for each pair of cell stabilisers,
using a presentation of the latter in terms of matrices, permutation cycles or generators and relators.
In the below examples however,
 we do avoid this case-by-case computation by a general determination of the isomorphism types of pairs of cell stabilisers
 for which group inclusion induces an isomorphism on mod $\ell$ cohomology.
The latter method is the procedure of preference,
because it allows us to deduce statements that hold for the entire class of groups in question.
\end{remark}

 \subsubsection{Example: A \texorpdfstring{$2$}{2}-torsion subcomplex for SL\texorpdfstring{$_3(\mathbb{Z})$}{(3,\textbf{Z})}}
 The $2$-torsion subcomplex of the cell complex described by Soul\'e~\cite{Soule},
 obtained from the action of SL$_3(\mathbb{Z})$ on its symmetric space,
 has the following homeomorphic image.
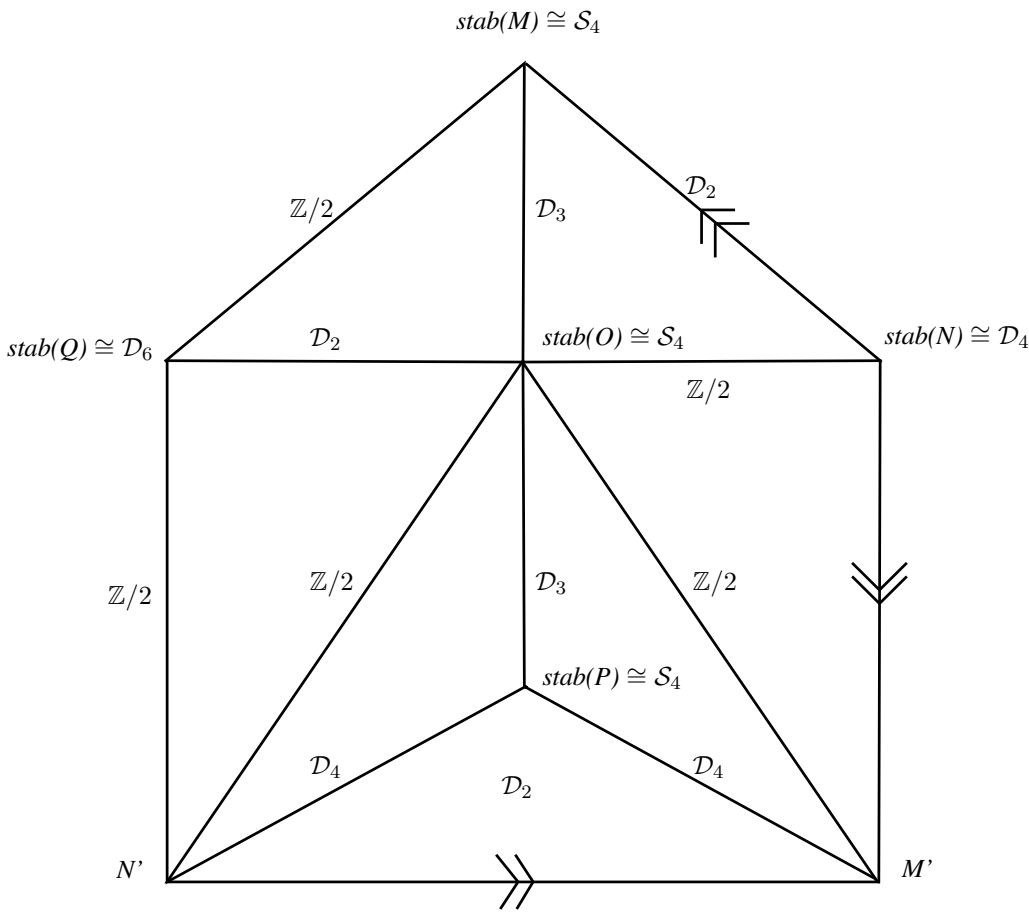
\begin{center}
\scalebox{0.9} 
{
\begin{pspicture}(-1.3,-7.44125)(11.894688,7.46125)
\pstriangle[linewidth=0.04,dimen=outer](5.8803124,-6.42125)(10.46,7.72)
\psline[linewidth=0.04](11.050312,-6.36125)(5.9103127,-3.52125)(5.8903127,1.25875)(0.6903125,1.27875)(0.6903125,-6.40125)(5.9303126,-3.52125)(5.9503126,-3.54125)
\psline[linewidth=0.04](11.090313,-6.42125)(11.110312,1.27875)(5.8903127,1.25875)(5.9103127,5.65875)(0.6703125,1.27875)(0.6903125,1.25875)
\psline[linewidth=0.04](5.9303126,5.65875)(11.110312,1.27875)(11.130313,1.25875)
\usefont{T1}{ptm}{m}{it}
\rput(5.9759374,6.26875){stab(M) $\cong \Sfour$}
\usefont{T1}{ptm}{m}{it}
\rput(-0.575,1.44875){stab(Q) $\cong \Dsix$}
\usefont{T1}{ptm}{m}{it}
\rput(7.217656,1.60875){stab(O) $\cong \Sfour$}
\usefont{T1}{ptm}{m}{it}
\uput[0](11.0,1.62875){stab(N) $\cong \Dfour$}
\usefont{T1}{ptm}{m}{it}
\uput[0](6.0,-3.4){stab(P) $\cong \Sfour$}
\usefont{T1}{ptm}{m}{it}
\rput(0.14765625,-6.21125){N'}
\usefont{T1}{ptm}{m}{it}
\rput(11.657657,-6.19125){M'}
\uput[90](8.5,3.5){$\Dtwo$} 
\uput[0](5.9,3.5){$\Dthree$} 
\uput[0](5.9,-2.0){$\Dthree$} 
\uput[0](5.4,-5.0){$\Dtwo$} 
\uput[180](3.3,3.5){$\Z/2$} 
\rput(0.14765625,-2.21125){$\Z/2$} 
\uput[0](2.6,-2.0){$\Z/2$} 
\uput[0](2.6,-4.7){$\Dfour$} 
\uput[0](8.2,-2.0){$\Z/2$} 
\uput[0](8.2,-4.7){$\Dfour$} 
\uput[0](2.6,1.6){$\Dtwo$} 
\uput[270](8.6,1.2){$\Z/2$} 
\psline[linewidth=0.04](8.99,3.5)(8.5,3.5)(8.5,3.0)
\psline[linewidth=0.04](9.2,3.3)(8.7,3.3)(8.7,2.8)
\psline[linewidth=0.04](10.7,-1.7)(11.1,-2.1)(11.5,-1.7)
\psline[linewidth=0.04](10.7,-1.9)(11.1,-2.3)(11.5,-1.9)
\psline[linewidth=0.04](5.49875,-6.0)(5.81875,-6.4)(5.55875,-6.8)
\psline[linewidth=0.04](5.77875,-6.0)(6.03875,-6.4)(5.79875,-6.8)
\end{pspicture}
}
\end{center}
Here, the three edges $NM$, $NM'$ and $N'M'$ have to be identified as indicated by the arrows.
All of the seven triangles belong with their interior to the $2$-torsion subcomplex,
each with stabiliser $\Z/2$, except for the one which is marked to have stabiliser $\Dtwo$.
Using the methods described in Section~\ref{conjugacy reduction},
we reduce this subcomplex to

\begin{center}
 \scalebox{1} 
{
\begin{pspicture}(-1.9,-0.9)(8.5,0.3)
        \psdots(-0.0,0.0)
        \psline(-0.0,0.0)(2.0,0.0)
        \uput{0.1}[90](-0.0,0.0){ $\Sfour$}
                \uput{0.4}[270](-0.1,0.2){ $O$}
        \psdots(2,0.0)
                \uput{0.1}[90](1.0,0.0){ $\Dtwo$}
        \uput{0.1}[90](2.0,0.0){ $\Dsix$}
                \uput{0.4}[270](2.0,0.2){ $Q$}
        \psline(2,0.0)(8,0.0)
                        \uput{0.1}[90](3.0,0.0){ $\Z/2$}
        \uput{0.1}[90](4.0,0.0){ $\Sfour$}
                \uput{0.4}[270](4.0,0.2){$M$}
        \psdots(4,0.0)
                        \uput{0.1}[90](5.0,0.0){ $\Dfour$}
               \uput{0.1}[90](6.0,0.0){ $\Sfour$}
                  \uput{0.2}[270](6.0,0.0){$P$}
                          \psdots(6,0.0)
                      \uput{0.1}[90](7.0,0.0){ $\Dfour$}
        \uput{0.1}[90](8.0,0.0){ $\Dfour$}
                \uput{0.4}[270](8.0,0.2){$N'$}
        \psdots(8,0.0)
\end{pspicture}
}
\end{center}
and then to
\begin{center}
 \scalebox{1} 
{
\begin{pspicture}(-1.9,-0.2)(8.0,0.3)
        \uput{0.1}[90](2.0,0.0){ $\Sfour$}
        \psdots(2,0.0)
        \psline(2,0.0)(6,0.0)
                        \uput{0.1}[90](3.0,0.0){ $\Z/2$}
        \uput{0.1}[90](4.0,0.0){ $\Sfour$}
        \psdots(4,0.0)
                        \uput{0.1}[90](5.0,0.0){ $\Dfour$}
               \uput{0.1}[90](6.0,0.0){ $\Sfour$}
        \psdots(6,0.0)
\end{pspicture}
}
\end{center}
which is the geometric realization of Soul\'e's diagram of cell stabilisers.
This yields the mod $2$ Farrell--Tate cohomology as specified in~\cite{Soule}.

\subsubsection{Example: Farrell--Tate cohomology of the Bianchi modular groups}
Consider the $\SLtwo$ matrix groups over the ring $\ringOm$
 of integers in the imaginary quadratic number field $\rationals(\sqrt{-m})$,
 with $m$ a square-free positive integer.
These groups, as well as their central quotients $\PSLO$, are known as \textit{Bianchi (modular) groups}.
We recall the following information from~\cites{Rahm:formulas} on the $\ell$-torsion subcomplex of $\PSLO$. Let~$\Gamma$ be a finite index subgroup in $\text{PSL}_2(\mathcal{O}_{-m})$.
Then any element of~$\Gamma$ fixing a point inside hyperbolic $3$-space~$\hyper$ acts as a rotation of finite order.
By Felix Klein's work, we know conversely that any torsion element~$\alpha$ is elliptic and hence fixes some geodesic line.
We call this line \emph{the rotation axis of~$\alpha$}.
Every torsion element acts as the stabiliser of a line conjugate to one passing through the Bianchi fundamental polyhedron.
We obtain the \textit{refined cellular complex} from the action of~$\Gamma$ on~$\hyper$
 as described in~\cite{Rahm:homological_torsion},
namely we subdivide~$\hyper$  until the stabiliser in~$\Gamma$ of any cell $\sigma$ fixes $\sigma$ point-wise.
We achieve this by computing Bianchi's fundamental polyhedron for the action of~$\Gamma$,
 taking as a preliminary set of 2-cells its facets lying on the Euclidean hemispheres
 and vertical planes of the upper-half space model for $\hyper$,
 and then subdividing along the rotation axes of the elements of~$\Gamma$.

It is well-known~\cite{SchwermerVogtmann} that if $\gamma$ is an element of finite order $n$ in a Bianchi group, then $n$ must be 1, 2, 3, 4 or 6,
 because $\gamma$ has eigenvalues $\rho$ and $\overline{\rho}$,
 with $\rho$ a primitive $n$-th root of unity, and the trace of~$\gamma$ is $\rho + \overline{\rho} \in \calO_{-m} \cap {\mathbb{R}} = \Z$.
 When $\ell$ is one of the two occurring prime numbers $2$ and~$3$, the orbit space of this subcomplex is a graph,
 because the cells of dimension greater \mbox{than 1} are trivially stabilized in the refined cellular complex.
We can see that this graph is finite either from the finiteness of the Bianchi fundamental polyhedron,
or from studying conjugacy classes of finite subgroups as in~\cite{Kraemer:Diplom}.

As in \cite{RahmFuchs}, we make use of a $2$-dimensional deformation retract $X$ of the refined cellular complex,
equivariant with respect to a Bianchi group \mbox{$\Gamma$}.  This retract has a cell structure
in which each cell stabiliser fixes its cell pointwise.
Since $X$ is a deformation retract of $\hyper$ and hence acyclic,
$$\Cohomol^*_\Gamma(X) \cong \Cohomol^*_\Gamma(\hyper) \cong \Cohomol^*(\Gamma).$$
\begin{table}
\begin{center}
$
\begin{array}{|c|c|c|c|c|c|}
\hline \text{Subgroup type}   & \Z/2 & \Z/3 & \Dtwo &\Dthree& \Af \\
\hline \text{Number of conjugacy classes} &  \lambda_{4} &  \lambda_6  &  \mu_2      & \mu_3   & \mu_T   \\
\hline
\end{array}
$
\end{center}
\caption{The non-trivial finite subgroups of $\PSLO$ were classified by Klein~\cite{Klein:binaereFormenMathAnn9}.
Here, $\Z/n$ is the cyclic group of order $n$, the dihedral groups are $\Dtwo$ with four elements and $\Dthree$ with six elements,
and the tetrahedral group is isomorphic to the alternating group $\Af$ on four letters.
Formulas for the numbers of conjugacy classes counted by the Greek symbols are given by Kr\"amer~\cite{Kraemer:Diplom}.}
\label{table:covering}
\end{table}
In Theorem~\ref{Grunewald-Poincare series formulas} below,
we give a formula expressing precisely how the Farrell--Tate cohomology
 of a Bianchi group with units $\{\pm 1\}$
 (i.e., just excluding the Gaussian and the Eisentein integers as imaginary quadratic rings,
 see Section~\ref{Bredon state})
 depends on the numbers of conjugacy classes of non-trivial finite subgroups of the occurring five types specified in Table~\ref{table:covering}.
The main step in order to prove this, is to read off the Farrell--Tate cohomology from the quotient of the reduced torsion sub-complexes.

Kr\"amer's formulas \cite{Kraemer:Diplom} express the numbers of conjugacy classes of the five types of non-trivial finite subgroups given in Table~\ref{table:covering}.
We are going to use the symbols of that table also for the numbers of conjugacy classes in $\Gamma$,
 where $\Gamma$ is a finite index subgroup in a Bianchi group.
Recall that for $\ell = 2$ and $\ell = 3$,
we can express the dimensions of the homology of $\Gamma$ with coefficients in the field ${\mathbb{F}_\ell}$ with $\ell$ elements
 in degrees above the virtual cohomological dimension of the Bianchi groups -- which is $2$ -- by the Poincar\'e series
$$P^\ell_\Gamma(t) := \sum\limits_{q \thinspace > \thinspace 2}^{\infty} \dim_{\mathbb{F}_\ell} \Homol_q \left(\Gamma;\thinspace {\mathbb{F}_\ell} \right)\thinspace t^q,$$
which has been suggested by Grunewald.
Further let $P_{\circlegraph} (t) := \frac{-2t^3}{t-1}$ , which equals the Poincar\'e series $P^2_\Gamma(t)$  of the groups $\Gamma$ the quotient of the reduced $2$--torsion sub-complex of which is a circle.
Denote by \begin{itemize}
\item $P_{\Dtwo}^*(t) := \frac{-t^3(3t -5)}{2(t-1)^2}$, the Poincar\'e series over
$$\dim_{\mathbb{F}_2} \Homol_q \left(\Dtwo;\thinspace {\mathbb{F}_2} \right) -\frac{3}{2}\dim_{\mathbb{F}_2} \Homol_q \left(\Z/2;\thinspace {\mathbb{F}_2} \right)$$
\item and by $P_{\Afour}^*(t) := \frac{-t^3(t^3 - 2t^2 + 2t - 3)}{2(t-1)^2 (t^2 + t + 1 ) }$, the Poincar\'e series over
$$\dim_{\mathbb{F}_2} \Homol_q \left(\Afour;\thinspace {\mathbb{F}_2} \right) -\frac{1}{2}\dim_{\mathbb{F}_2} \Homol_q \left(\Z/2;\thinspace {\mathbb{F}_2} \right).$$
\end{itemize}

In 3-torsion, let
$P_{\edgegraph} (t) := \frac{-t^3(t^2 - t + 2)}{(t-1)(t^2+1)}$, which equals the Poincar\'e series $P^3_\Gamma(t)$ for those Bianchi groups,
the quotient of the reduced $3$--torsion sub-complex of which is a single edge without identifications.

\vbox{
\begin{theorem}[\cite{Rahm:formulas}] \label{Grunewald-Poincare series formulas}
For any finite index subgroup $\Gamma$ in a Bianchi group with units $\{\pm 1\}$, the group homology in degrees above its virtual cohomological dimension is given by the Poincar\'e series
$$P^2_\Gamma(t) = \left(\lambda_4 -\frac{3\mu_2 -2\mu_T}{2}\right)P_{\circlegraph} (t) +(\mu_2 -\mu_T)P_{\Dtwo}^*(t) +\mu_T P_{\Afour}^*(t)$$
and
$$P^3_\Gamma(t) =  \left(\lambda_6 -\frac{\mu_3}{2}\right)P_{\circlegraph} (t) + \frac{\mu_3}{2}P_{\edgegraph}(t).$$
\end{theorem}
}

More general results are stated in Section~\ref{formulas for the Farrell--Tate cohomology} below.

 \subsubsection{Example: Farrell--Tate cohomology of Coxeter (tetrahedral) groups}
Recall that a Coxeter group is a group admitting a presentation
$$\langle g_1, g_2, ..., g_n \medspace | \medspace (g_i g_j)^{m_{i,j}} = 1 \rangle,$$
where $m_{i,i} = 1$; for $i \neq j$ we have $m_{i,j} \geq 2$;
 and $m_{i,j} = \infty$ is permitted, meaning that $(g_i g_j)$ is not of finite order.
As the Coxeter groups admit a contractible classifying space for proper actions \cite{Davis},
 their Farrell--Tate cohomology yields all of their group cohomology.
So in this section, we make use of this fact to determine the latter.
For facts about Coxeter groups, and especially for the Davis complex, we refer to \cite{Davis}.
Recall that the simplest example of a Coxeter group, the dihedral group $\Dn$, is an extension
$$ 1 \to \Z/n \to \Dn \to \Z/2 \to 1.$$
So we can make use of the original application~\cite{Wall} of Wall's lemma to obtain its mod $\ell$ homology for prime numbers
 $\ell >2$,
$$ \Homol_q(\Dn; \thinspace \Z/\ell) \cong
\begin{cases}
				\Z/\ell, & q = 0, \\
                               \Z/{\rm gcd}(n,\ell), & q \equiv 3 \medspace {\rm or} \medspace 4 \mod 4, \\
				0, & {\rm otherwise}.
\end{cases}
$$
\begin{theorem}[\cite{Rahm:formulas}] \label{small rank Coxeter groups}
 Let $\ell > 2$ be a prime number.
 Let $\Gamma$ be a Coxeter group admitting a Coxeter system with at most four generators,
 and relator orders not divisible by~$\ell^2$.
Let $Z_{(\ell)}$ be the $\ell$--torsion sub-complex of the Davis complex of~$\Gamma$.
If $Z_{(\ell)}$ is at most one-dimensional and its orbit space contains no loop or bifurcation, then the$\mod \ell$ homology of~$\Gamma$ is isomorphic to
$\left(\Homol_q(\Dl; \thinspace \Z/\ell)\right)^m$, with $m$ the number of connected components of the orbit space of~$Z_{(\ell)}$.
\end{theorem}
The conditions of this theorem are for instance fulfilled by the Coxeter tetrahedral groups;
 the exponent $m$ has been specified for each of them in the tables in~\cite{Rahm:formulas}.
In the easier case of Coxeter triangle groups, we can sharpen the statement as follows.
 The non-spherical and hence infinite \emph{Coxeter triangle groups} are given by the presentation
$$
    \langle\, a, b, c \;|\; a^2 = b^2 = c^2 = (ab)^p = (bc)^q =
    (c a)^r = 1 \,\rangle\, ,
$$
where $2 \leq p,q,r \in \N$ and $\frac{1}{p} + \frac{1}{q} +
\frac{1}{r} \le 1$.

\begin{proposition}[\cite{Rahm:formulas}]
For any prime number $\ell>2$, the {\rm mod} $\ell$ homology of a Coxeter triangle group is given as the direct sum over
the {\rm mod} $\ell$ homology of the dihedral groups ${\mathcal D}_p$, ${\mathcal D}_q$ and ${\mathcal D}_r$.
\end{proposition}

\subsection{The non-central torsion subcomplex} \label{The non-central torsion subcomplex}

In the case of a trivial kernel of the action on the polytopal $\Gamma$-cell complex,
torsion subcomplex reduction allows one to establish general formulas for the Farrell--Tate cohomology of~$\Gamma$ \cite{Rahm:formulas}.
In contrast, for instance the action of $\SLO$ on hyperbolic $3$-space has the $2$-torsion group $\{\pm 1\}$ in the kernel;
since  every cell stabiliser contains  $2$-torsion,
the $2$-torsion subcomplex then does not ease our calculation in any way.
We can remedy this situation by considering the following object, on whose cells we impose a supplementary property.

\begin{df} \label{non-central torsion subcomplex}
 The \emph{non-central $\ell$-torsion subcomplex} of a polytopal $\Gamma$-cell complex $X$
 consists of all the cells of $X$
 whose stabilisers in~$\Gamma$  contain elements of order $\ell$ that are not in the center of~$\Gamma$.
\end{df}

We note that this definition yields a correspondence between, on one side, the \textit{non-central}
$\ell$-torsion subcomplex for a group action with kernel the center of the group,
 and on the other side, the $\ell$-torsion subcomplex for its central quotient group.
In~\cite{BerkoveRahm}, this correspondence has been used in order to identify the \textit{non-central} $\ell$-torsion subcomplex for the action of $\SLO$
 on hyperbolic $3$-space as the $\ell$-torsion subcomplex of $\PSLO$.
  However, incorporating the non-central condition for $\SLO$
  introduces significant technical obstacles,
  which were addressed in that paper, establishing
  the following theorem for any finite index subgroup $\Gamma$ in $\SLO$.
Denote by $X$ a $\Gamma$-equivariant retract of SL$_2(\C)/$SU$_2$,
by $X_s$ the $2$-torsion subcomplex with respect to P$\Gamma$
(the ``non-central'' $2$-torsion subcomplex for $\Gamma$),
and by $\xsp$ the part of it with higher $2$-rank.
Further, let $v$ denote the number of conjugacy classes of subgroups of higher $2$-rank, and define
$\sign(v) := $\scriptsize$\begin{cases}
                   0, & v = 0,\\
		    1,& v> 0.
                  \end{cases}$\normalsize
                  \\
For $q \in \{1, 2\}$, denote the dimension $\dim_{\F_2}\Cohomol^q(_\Gamma \backslash X ; \thinspace \F_2)$ by $\beta^q$.

\vbox{ \begin{thm}[\cite{BerkoveRahm}] \label{E2 page}
 The $E_2$ page of the equivariant spectral sequence with $\F_2$-coefficients
 associated to the action of $\Gamma$ on $X$ is concentrated in the columns $n \in \{0, 1, 2\}$
and has the following form.
 \[
\begin{array}{l | cccl}
q = 4k+3 &  E_2^{0,3}(X_s)     &   E_2^{1,3}(X_s) \oplus (\ef)^{a_1}  &   (\ef)^{a_2} \\
q = 4k+2 &  \Cohomol^2_\Gamma(\xsp) \oplus (\ef)^{1 -\sign(v)}     &    (\ef)^{a_3}&   \Cohomol^2(_\Gamma \backslash X)  \\
q = 4k+1 &  E_2^{0,1}(X_s)     &  E_2^{1,1}(X_s) \oplus  (\ef)^{a_1} &   (\ef)^{a_2} \\
q = 4k   & \F_2  &   \Cohomol^1(_\Gamma \backslash X) &  \Cohomol^2(_\Gamma \backslash X) \\
\hline k \in \mathbb{N} \cup \{0\} &  n = 0 & n = 1 &  n = 2
\end{array}
\]
where
\[\begin{array}{ll}
a_1 & =  \chi(_\Gamma \backslash X_s) -1 +\beta^1(_\Gamma \backslash X) +c   \\
a_2 & =  \beta^{2} (_\Gamma \backslash X) +c \\
a_3 & =  \beta^{1} (_\Gamma \backslash X) +v -\sign(v).
\end{array}
\]
\end{thm}
}

In order to derive the example stated in Section~\ref{non-trivial-centre} below,
we combine the latter theorem with the following determination (carried out in~\cite{BerkoveRahm})
of the $d_2$-differentials on the four possible (cf. Table~\ref{table:subcomplexes}) connected component types
$\circlegraph$, $\edgegraph$, $\graphFive$ and $\graphTwo$
of the reduced non-central $2$-torsion subcomplex for the full SL$_2$
groups over the imaginary quadratic number rings.
\begin{lemma}[\cite{BerkoveRahm}] \label{d_2 lemma}
The $d_2$ differential in the equivariant spectral sequence associated to the action of
$\SLtwo(\ringOm)$ on hyperbolic space is trivial on components of the non-central
$2$-torsion subcomplex quotient
\begin{itemize}
 \item of type $\circlegraph$ in dimensions $q \equiv 1 \bmod 4$
if and only if it is trivial on these components in dimensions $q \equiv 3 \bmod 4$.
\item of type $\edgegraph$.
\item of types $\graphTwo$ and $\graphFive$ in dimensions $q \equiv 3 \bmod 4$.
\end{itemize}
\end{lemma}

\begin{table}
 \begin{center}
 \caption{Connected component types of reduced torsion subcomplex quotients for the PSL$_2$ Bianchi groups.
 The exhaustiveness of this table has been established using theorems of Kr\"amer \cite{BerkoveRahm}.} \label{table:subcomplexes}
 \label{one}
 \footnotesize
 \begin{tabular}{|c|c|c|c|c|}
 \hline & & &&\\
\begin{tabular}{c}$2$--torsion\\subcomplex \\components\end{tabular}
& \begin{tabular}{c}counted \\by \end{tabular} & &
\begin{tabular}{c}$3$--torsion\\subcomplex \\components\end{tabular}
& \begin{tabular}{c}counted \\by \end{tabular}
\\  \hline & & &&\\
 $\circlegraph \thinspace \Z/2$ & $o_2 = \lambda_4 -\lambda_4^* $ & &  $\circlegraph \thinspace \Z/3$ & $o_3  = \lambda_6 -\lambda_6^* $\\
 & & & &\\
$\Afour \edgegraph \Afour$ & $\iota_2$ & & $\Dthree \edgegraph \Dthree$ & $\iota_3 = \lambda_6^* $\\
 & & &&\\
$\Dtwo \graphFive \thinspace \Dtwo$ & $\theta$&&&\\
 & & &&\\
$\Dtwo \graphTwo \Afour$ & $\rho$ &&&\\
\hline
\end{tabular}
\normalsize
\end{center}
\end{table}

\section{Applications of the technique and their results} \label{Results}
This section is going to state some results in which the technique described in Section~\ref{techniques} was involved.

\subsection{The Bianchi groups and their congruence subgroups} \label{The Bianchi groups}
In the case of the PSL$_2$ groups over rings of imaginary quadratic integers (known as the Bianchi groups),
the torsion subcomplex reduction technique
has permitted the author to find a description of the cohomology ring of these groups in terms of elementary number-theoretic quantities~\cite{Rahm:formulas}.
The key step has been to extract, using torsion subcomplex reduction, the essential information about the geometric models, and then to detach the cohomological information completely from the model.

Torsion subcomplex reduction combined with an analysis of the equivariant spectral sequence by Ethan Berkove, Grant Lakeland and the author provides new tools for the calculation of the torsion in the cohomology of congruence subgroups in the Bianchi groups~\cite{BLR}.

\subsection{The Coxeter groups} \label{The Coxeter groups}
Let us recall that the Coxeter groups are generated by reflections,
and their homology consists solely of torsion.
Thus, torsion subcomplex reduction allows one to obtain all homology groups for all of the tetrahedral Coxeter groups
at all odd prime numbers, in terms of a general formula~\cite{Rahm:formulas}.

\subsection{The SL\texorpdfstring{$_2$}{(2)} groups over arbitrary number rings}\label{formulas for the Farrell--Tate cohomology}
Matthias Wendt and the author established a complete
description of the Farrell--Tate cohomology with odd torsion coefficients
for all groups $\op{SL}_2(\mathcal{O}_{K,S})$,
where $\mathcal{O}_{K,S}$ is the ring of $S$-integers in an arbitrary number field $K$
at an arbitrary non-empty finite set $S$ of places of $K$ containing the infinite places~\cite{RahmWendt},
based on an explicit description of
conjugacy classes of finite cyclic subgroups and their normalizers in
$\op{SL}_2(\mathcal{O}_{K,S})$.

The statement uses the following notation.
Let $\ell$ be an odd prime number different from the characteristic of $K$.
In the situation where, for $\zeta_\ell$ some primitive $\ell$-th root of unity,
$\zeta_\ell+\zeta_\ell^{-1}\in K$,
we will abuse notation and write $\mathcal{O}_{K,S}[\zeta_\ell]$ to mean the ring
$\mathcal{O}_{K,S}[T]/(T^2-(\zeta_\ell+\zeta_\ell^{-1})T+1)$.
Moreover, we denote the norm maps for class groups and units by
$$
\op{Nm}_0:
  \widetilde{\op{K}_0}(\mathcal{O}_{K,S}[\zeta_\ell])\to
  \widetilde{\op{K}_0}(\mathcal{O}_{K,S})
\qquad\textrm{ and }\qquad
\op{Nm}_1:\mathcal{O}_{K,S}[\zeta_\ell]^\times\to
  \mathcal{O}_{K,S}^\times.
$$
Denote by $M_{(\ell)}$ the $\ell$-primary part of a module $M$;
by $N_G(\Gamma)$ the normalizer of $\Gamma$ in $G$;
and by $\widehat{\op{H}}^\bullet$ Farrell--Tate cohomology (cf. Section~\ref{background}).

\begin{theorem} [\cite{RahmWendt}]
\label{thm:gl2nf}
${}$
\begin{enumerate}
\item
  $\widehat{\op{H}}^\bullet(\op{SL}_2(\mathcal{O}_{K,S}),\mathbb{F}_\ell)\neq
  0$ if and only if \\
  $\zeta_\ell+\zeta_\ell^{-1}\in K$ and the Steinitz
  class $\det_{\mathcal{O}_{K,S}}(\mathcal{O}_{K,S}[\zeta_\ell])$ is
  contained in the image of the norm map $\op{Nm}_0$.
\item Assume the condition in (1) is satisfied.
  The set $\mathcal{C}_\ell$ of conjugacy classes of order $\ell$
  elements in $\op{SL}_2(\mathcal{O}_{K,S})$
  sits in an extension
$$
1\to \op{coker}\op{Nm}_1\to \mathcal{C}_\ell\to
\ker\op{Nm}_0\to 0.
$$
The set $\mathcal{K}_\ell$ of conjugacy classes of order $\ell$
subgroups of $\op{SL}_2(\mathcal{O}_{K,S})$ can be identified with the
quotient
$\mathcal{K}_\ell=\mathcal{C}_\ell/\op{Gal}(K(\zeta_\ell)/K)$. There is
a direct sum decomposition
$$
\widehat{\op{H}}^\bullet(\op{SL}_2(\mathcal{O}_{K,S}),\mathbb{F}_\ell)\cong
\bigoplus_{[\Gamma]\in\mathcal{K}_\ell}
\widehat{\op{H}}^\bullet(N_{\op{SL}_2(\mathcal{O}_{K,S})}(\Gamma),
\mathbb{F}_\ell)
$$
which is compatible with the ring structure, i.e., the Farrell--Tate cohomology ring of $\op{SL}_2(\mathcal{O}_{K,S})$ is a direct sum of the sub-rings
for the normalizers $N_{\op{SL}_2(\mathcal{O}_{K,S})}(\Gamma)$.

\item If the class of $\Gamma$ is not $\op{Gal}(K(\zeta_\ell)/K)$-invariant,
then
$$N_{\op{SL}_2(\mathcal{O}_{K,S})}(\Gamma)\cong \ker\op{Nm}_1.$$
There is a degree $2$ cohomology class $a_2$ and a ring isomorphism
$$
\widehat{\op{H}}^\bullet(\ker\op{Nm}_1,\mathbb{Z})_{(\ell)}\cong
\mathbb{F}_\ell[a_2,a_2^{-1}]\otimes_{\mathbb{F}_\ell}\bigwedge
\left(\ker\op{Nm}_1\right).
$$ In particular, this is a free module
over the subring $\mathbb{F}_\ell[a_2^2,a_2^{-2}]$.
\item
If the class of $\Gamma$ is $\op{Gal}(K(\zeta_\ell)/K)$-invariant,
then there is an extension
$$
0\to \ker\op{Nm}_1\to N_{\op{SL}_2(\mathcal{O}_{K,S})}(\Gamma)\to \mathbb{Z}/2\to 1.
$$
There is a ring isomorphism
$$
\widehat{\op{H}}^\bullet(N_{\op{SL}_2(\mathcal{O}_{K,S})}(\Gamma),\mathbb{Z})_{(\ell)}\cong
\left(\mathbb{F}_\ell[a_2,a_2^{-1}]\otimes_{\mathbb{F}_\ell}\bigwedge
\left(\ker\op{Nm}_1\right)\right)^{\mathbb{Z}/2},
$$
with the $\mathbb{Z}/2$-action given by multiplication with $-1$ on
$a_2$ and $\ker\op{Nm}_1$. In particular, this is a free module over
the subring
$$\mathbb{F}_\ell[a_2^2,a_2^{-2}]\cong
\widehat{\op{H}}^\bullet(D_{2\ell},\mathbb{Z})_{(\ell)}.$$
\item The restriction map induced from the inclusion
  $\op{SL}_2(\mathcal{O}_{K,S})\to
  \op{SL}_2(\mathbb{C})$ maps the second Chern class
  $c_2$ to the sum of the elements $a_2^2$ in all the components.
\end{enumerate}
\end{theorem}
Wendt has extended this investigation to the cases of $\op{SL}_2$ over the ring of functions on a smooth affine curve over an algebraically closed field~\cite{sl2parabolic}.

\subsection{Farrell--Tate cohomology of higher rank arithmetic groups} \label{GL3}
Pertinent progress was also made on the Farrell--Tate cohomology of GL$_3$ over rings of quadratic integers.
For this purpose, the conjugacy classification of cyclic subgroups was reduced to the classification of modules of group rings over suitable rings of integers which are principal ideal domains, generalizing an old result of Reiner. As an example of the number-theoretic input required for the Farrell--Tate cohomology computations, Bui, Wendt and the author describe the homological torsion in PGL$_3$ over principal ideal rings of quadratic integers, accompanied by machine computations in the imaginary quadratic case~\cite{BuiRahmWendt:GL3om}.

For machine calculations of Farrell--Tate or Bredon (co)homology, one needs cell complexes where cell stabilizers fix their cells pointwise.
Bui, Wendt and the author provided two algorithms computing an efficient subdivision of a complex to achieve this rigidity property~\cite{BuiRahmWendt:Farrell-Tate}.
Applying these algorithms to available cell complexes for {PSL}$_4(\Z)$,
they computed the Farrell--Tate cohomology for small primes as well as the Bredon homology for the classifying spaces of proper actions with coefficients in the complex representation ring.

\subsection{Adaptation of the technique to groups with non-trivial centre}
\label{non-trivial-centre}
Berkove and the author~\cite{BerkoveRahm} extended the technique of torsion subcomplex reduction,
which originally was designed for groups with trivial centre (e.g., PSL$_2$),
to groups with non-trivial centre (e.g., SL$_2$).
This way, they determined the $2$-torsion in the cohomology of the SL$_2$ groups
over imaginary quadratic number rings~$\ringOm$ in $\rationals(\sqrt{-m})$,
based on their action on hyperbolic 3-space~$\Hy$.
For instance, they get the following result in the case where the quotient of the
$2$--torsion subcomplex has the shape $\edgegraph$,
which is equivalent to the following three conditions (cf.~\cite{Rahm:formulas}):
$m \equiv 3 \bmod 8$, the field $\rationals(\sqrt{-m})$ has precisely one finite ramification place over $\rationals$,
and the ideal class number of the totally real number field $\rationals(\sqrt{m})$ is $1$.
Under these assumptions, our cohomology ring has the following dimensions:
$$
\dim_{\F_2}\Cohomol^{q}(\SLO; \thinspace \F_2)
=
\begin{cases}
   \beta^{1}  +\beta^{2} , &  q = 4k+5, \\
   \beta^1 +\beta^{2} +2, &  q = 4k+4, \\
   \beta^{1} + \beta^{2}+3, &  q = 4k+3, \\
   \beta^1+\beta^{2} +1, &  q = 4k+2, \\
   \beta^{1},  &  q = 1, \\
\end{cases}
$$
where $\beta^q := \dim_{\F_2}\Cohomol^q(_{\SLO} \backslash \Hy ; \thinspace \F_2)$.
Let $ \beta_1 := \dim_{\rationals}\Cohomol_1(_{\SLO} \backslash \Hy ; \thinspace \rationals).$
For all absolute values of the discriminant less than $296$, numerical calculations yield
$\beta^2 +1 = \beta^1= \beta_1.$
In this range, the numbers $m$ subject to the above dimension formula and $\beta_1$ are given as follows
 (the Betti numbers are computed in a previous paper of the author~\cite{Rahm:higher_torsion}).
$$\begin{array}{l|ccccccccccccccccccccccccccc}
m        &  11 & 19 & 43 & 59 & 67 & 83 & 107 & 131 & 139 & 163 & 179 & 211 & 227 & 251 & 283\\
 \hline
\beta_1  &   1 & 1  &  2 &  4 & 3  &  5 &  6  &   8 &  7  &  7  &  10 &  10 & 12  & 14  & 13\\
\end{array}$$

This result is a consequence of Theorem~\ref{E2 page},
combined with Lemma~\ref{d_2 lemma} above.

\subsection{Investigation of the refined Quillen conjecture} \label{QC}
The Quillen conjecture on the cohomology of arithmetic groups has spurred a great deal of mathematics
(see the pertinent monograph \cite{Knudson:book}).
Using Farrell--Tate cohomology computations, Wendt and the author established further positive cases for the Quillen
conjecture for $\op{SL}_2$.
In detail, the original conjecture of 1971~\cite{Quillen} is as follows for GL$_n$.
\begin{conjecture}[Quillen] \label{Quillen-conjecture}
Let $\ell$ be a prime number. Let $K$ be a number field with
$\zeta_\ell\in K$, and $S$ a finite set of places containing the
infinite places and the places over $\ell$. Then the natural inclusion
$\mathcal{O}_{K,S}\hookrightarrow \mathbb{C}$ makes
$\op{H}^\bullet(\op{GL}_n(\mathcal{O}_{K,S}),\mathbb{F}_\ell)$ a free
module over the cohomology ring
$\op{H}^\bullet_{\op{cts}}(\op{GL}_n(\mathbb{C}),\mathbb{F}_\ell)$.
\end{conjecture}
While there are counterexamples to the original version of the conjecture,
it holds true in many other cases.
From the first counterexamples through the present,
the conjecture has kept researchers interested in determining its range of validity~\cite{Anton-mod5}.

Positive cases in which the conjecture  has been established
are $n=\ell=2$ by Mitchell \cite{Mitchell}, $n=3$, $\ell=2$ by Henn
\cite{Henn}, and  $n=2$, $\ell=3$ by Anton \cite{anton}.

On the other hand, cases where the Quillen conjecture is known to be
false can all be traced to a remark by Henn, Lannes and Schwartz~\cite{henn:lannes:schwartz}*{remark on p. 51}, which shows that Quillen's conjecture
for $\op{GL}_n(\mathbb{Z}[1/2])$ implies that the restriction map
$$
\op{H}^\bullet(\op{GL}_n(\mathbb{Z}[1/2]),\mathbb{F}_2)\to
\op{H}^\bullet(\op{T}_n(\mathbb{Z}[1/2]),\mathbb{F}_2)
$$
from $\op{GL}_n(\mathbb{Z}[1/2])$ to the  subgroup
$\op{T}_n(\mathbb{Z}[1/2])$ of diagonal matrices is
injective. Non-injectivity of the restriction map has been shown by
Dwyer \cite{dwyer}  for $n\geq 32$ and $\ell=2$. Dwyer's bound
was subsequently improved by Henn and Lannes to $n\geq 14$. At the prime
$\ell=3$, Anton~\cite{anton} proved non-injectivity for $n\geq 27$.

Wendt's and the author's  contribution is that we can determine precisely the module
structure above the virtual cohomological dimension; this has allowed us to
relate the Quillen conjecture for $\op{SL}_2$  to statements about
Steinberg homology (recall Section~\ref{Farrell--Tate cohomology and Steinberg homology}).

This, together with the results of~\cite{sl2parabolic},
has allowed us to find a refined version of the Quillen conjecture,
which keeps track of all the types of known counter-examples to the original Quillen conjecture:

\begin{conjecture}[Refined Quillen conjecture \cite{qcnote}]
\label{refined Quillen-conjecture}
Let $K$ be a number field. Fix a prime $\ell$  such that
$\zeta_\ell\in K$, and an integer $n<\ell$. Assume that $S$ is a set
of places containing the infinite places and the places lying over
$\ell$.
If each cohomology class of $\op{GL}_n(\mathcal{O}_{K,S})$ is
detected on some finite subgroup, then
$\op{H}^\bullet(\op{GL}_n(\mathcal{O}_{K,S}),\mathbb{F}_\ell)$ is a
free module over the image of the restriction map
$$\op{H}^\bullet_{\op{cts}}(\op{GL}_n(\mathbb{C}),\mathbb{F}_\ell)\to
\op{H}^\bullet(\op{GL}_n(\mathcal{O}_{K,S}),\mathbb{F}_\ell).$$
\end{conjecture}

We can make the following use of the description of the Farrell--Tate cohomology of SL$_2$ over rings of $S$-integers.

\begin{corollary}[Corollary to Theorem \ref{thm:gl2nf}]
Let $K$ be a number field, let $S$ be a  finite set of places
containing the infinite ones, and let $\ell$ be  an odd prime.
\begin{enumerate}
\item The original Quillen conjecture holds  for group cohomology
  $\op{H}^\bullet(\op{SL}_2(\mathcal{O}_{K,S}),\mathbb{F}_\ell)$ above
  the virtual cohomological dimension.
\item The refined Quillen conjecture holds for Farrell--Tate cohomology
  $\widehat{\op{H}}^\bullet(\op{SL}_2(\mathcal{O}_{K,S}),\mathbb{F}_\ell)$.
\end{enumerate}
\end{corollary}

\subsection*{Verification of the Quillen conjecture in the rank 2 imaginary quadratic case}
Bui and the author confirm a conjecture of Quillen in the case of the mod $2$ cohomology of arithmetic groups
${\rm SL}_2(\imQuadRing[\frac{1}{2}])$,
where $\imQuadRing$ is an imaginary quadratic ring of integers.
To make explicit the free module structure on the cohomology ring conjectured by Quillen, they computed the mod $2$ cohomology of
$\arithGrp$ via the amalgamated decomposition of the latter group~\cite{BuiRahm:Verification}.

\subsection{Application to equivariant \textit{K}-homology} \label{Bredon state}
For the Bianchi groups, the torsion subcomplex reduction technique was adapted from group homology to Bredon homology $\Homol^\mathfrak{Fin}_n(\Gamma; \thinspace R_\C)$
with coefficients in the complex representation rings, and with respect to the family of finite subgroups~\cite{Rahm:equivariant}.
This has led the author to the following formulas for this Bredon homology, and by the Atiyah--Hirzebruch spectral sequence,
to the below formulas for equivariant $K$-homology of the Bianchi groups acting on their classifying space for proper actions.

We let a Bianchi group $\Gamma$ act on a $2$-dimensional retract $X$ of hyperbolic $3$-space.
Denote by $\Gamma_\sigma$ the stabiliser of a cell $\sigma$, and by $R_\C(G)$ the complex representation ring of a group $G$.
As $X$ is a model for the classifying space for proper $\Gamma$-actions,
the homology of the Bredon chain complex of our $\Gamma$-cell complex $X$ is identical to the Bredon homology $\Homol^\mathfrak{Fin}_p(\Gamma; \thinspace R_\C)$
of~$\Gamma$~\cite{Sanchez-Garcia}. This Bredon chain complex can be stated as
$$\xymatrix{
0 \to \bigoplus\limits_{\sigma \in \thinspace_\Gamma \backslash X^{(2)}} R_\C (\Gamma_\sigma) \ar[r]^{\Psi_2 } &
\bigoplus\limits_{\sigma \in \thinspace_\Gamma \backslash X^{(1)}} R_\C (\Gamma_\sigma)  \ar[r]^{\Psi_1} &
\bigoplus\limits_{\sigma \in \thinspace_\Gamma \backslash X^{(0)}} R_\C (\Gamma_\sigma) \to 0,
}  $$
where the blocks of the differential matrices $\Psi_1$ and $\Psi_2$ are obtained by inducing homomorphisms on the involved complex representation rings from the cell stabiliser inclusions.
Note that in general, the Bredon chain complex continues in higher dimensions, and its truncation at dimension $2$ results from the dimension of $X$.

\begin{theorem} \label{splitting}
Let $\Gamma$ be a Bianchi group or any one of its subgroups.
 Then the Bredon homology $\Homol^\mathfrak{Fin}_n(\Gamma; \thinspace R_\C)$ splits as a direct sum over (1) the orbit space homology $\Homol_n(\underbar{\rm B}\Gamma; \thinspace \Z)$,
\begin{enumerate}
 \item[(2)] a submodule $\Homol_n(\Psi_\bullet^{(2)})$ determined by the reduced $2$-torsion subcomplex of $(\underline{\rm E}\Gamma, \Gamma)$ and
 \item[(3)] a submodule $\Homol_n(\Psi_\bullet^{(3)})$ determined by the reduced $3$-torsion subcomplex of $(\underline{\rm E}\Gamma, \Gamma)$.
\end{enumerate}
\end{theorem}
These submodules are given as follows.

Except for the Gauss{}ian and Eisenstein integers, which can easily be treated ad hoc~\cite{Rahm:noteAuxCRAS}, all the rings of integers of imaginary quadratic number fields admit as only units $\{\pm 1\}$. In the latter case, we call $\PSL_2(\ringOm)$ a \textit{Bianchi group with units} $\{\pm 1\}$.
\begin{theorem} \label{2}
 The $2$-torsion part of the Bredon complex of a {Bianchi group $\Gamma$ with units} $\{\pm 1\}$ has homology
 \begin{center}
  $\Homol_n(\Psi_\bullet^{(2)}) \cong \begin{cases}
                                      \Z^{z_2}\oplus (\Z/2)^\frac{d_2}{2},& n = 0,\\
                                      \Z^{o_2},& n = 1,\\
                                      0,&\text{\rm otherwise},
                                     \end{cases}
  $
 \end{center}
where $z_2$ counts the number of conjugacy classes of subgroups of type $\Z/2$ in $\Gamma$,
$o_2$ counts the conjugacy classes of type $\Z/2$ in $\Gamma$ which are not contained in any $2$-dihedral subgroup,
and $d_2$ counts the number of $2$-dihedral subgroups, whether or not they are contained in a tetrahedral subgroup of $\Gamma$.
\end{theorem}

\vbox{
\begin{theorem} \label{3}
 The $3$-torsion part of the Bredon complex of a {Bianchi group $\Gamma$ with units} $\{\pm 1\}$ has homology
 \begin{center}
  $\Homol_n(\Psi_\bullet^{(3)}) \cong \begin{cases}
                                      \Z^{2o_3+\iota_3},& n = 0 \medspace \text{\rm or }1,\\
                                      0,&\text{\rm otherwise},
                                     \end{cases}
  $
 \end{center}
where amongst the subgroups of type $\Z/3$ in $\Gamma$,
$o_3$ counts the number of conjugacy classes of those of them which are not contained in any $3$-dihedral subgroup,
and $\iota_3$ counts the conjugacy classes of those of them which are contained in some $3$-dihedral subgroup in $\Gamma$.
\end{theorem}
}

There are formulas for $o_2, z_2, d_2, o_3$ and $\iota_3$ in terms of elementary number-theoretic quantities~\cite{Kraemer:Diplom},
which are readily computable by machine~\cite{Rahm:formulas}*{appendix}. See Table~\ref{table:subcomplexes} for how they relate to the types of connected components of torsion subcomplexes.

We deduce in the following corollary formulas for the equivariant $K$-homology of the Bianchi groups. Note for this purpose that for a Bianchi group $\Gamma$, there is a model for \underline{E}$\Gamma$ of dimension 2, so $\Homol_2(\underline{\rm B}\Gamma ; \thinspace \Z) \cong \Z^{\beta_2}$ is torsion-free.
Note also that the naive Euler characteristic of the Bianchi groups vanishes (again excluding the two special cases of Gaussian and Eisensteinian integers), that is,
for $\beta_i = \dim \Homol_i(\underline{\rm B}\Gamma ; \thinspace \rationals)$ we have
$\beta_0 -\beta_1 +\beta_2 = 0$  and $\beta_0 = 1$.

Whenever we have a classifying space for proper $G$-actions of dimension at most $2$,
the Atiyah---Hirzebruch spectral sequence from its Bredon homology to the equivariant \mbox{$K$-homology} $K^G_j(\underbar{\rm E}G)$ of a group $G$
degenerates on the $E^2$-page and directly yields the following theorem, which can be found in the book by Mislin and Valette.
Note that by Bott periodicity, only two indices $j = 0$ and $j = 1$ are relevant.
\begin{thm}[\cite{MislinValette}] \label{Bredon_to_K-homology}
 Let $G$ be an arbitrary group such that $\dim \underbar{\rm E}G \leq 2$.
Then there is a natural short exact sequence
$$0 \to \Homol^\mathfrak{Fin}_0(G; R_\C) \to K^G_0(\underbar{\rm E}G) \to \Homol^\mathfrak{Fin}_2(G; R_\C) \to 0 $$
and a natural isomorphism $\Homol^\mathfrak{Fin}_1(G; R_\C) \cong K^G_1(\underbar{\rm E}G)$.
\end{thm}

\begin{corollary}[Corollary to theorems \ref{splitting}, \ref{2}, \ref{3} and \ref{Bredon_to_K-homology}]
 For any {Bianchi group $\Gamma$ with units} $\{\pm 1\}$, the short exact sequence linking Bredon homology and equivariant
 $K$-homology splits into
$$K^\Gamma_0(\underbar{\rm E}\Gamma) \cong \Z \oplus \Z^{\beta_2}  \oplus \Z^{z_2} \oplus (\Z/2)^\frac{d_2}{2} \oplus \Z^{2o_3+\iota_3}.$$
Furthermore, $K^\Gamma_1(\underbar{\rm E}\Gamma) \cong \Homol_1(\underline{\rm B}\Gamma; \thinspace \Z) \oplus \Z^{o_2} \oplus \Z^{2o_3+\iota_3}$.
\end{corollary}

In order to adapt torsion subcomplex reduction to Bredon homology and prove Theorem~\ref{splitting},
we need to perform a ``representation ring splitting''.

 \textit{Representation ring splitting}. \label{Representation ring splitting}
The classification of Felix Klein~\cite{Klein:binaereFormenMathAnn9}
of the finite subgroups in $\mathrm{PSL}_2(\ringO)$
 is recalled in Table~\ref{table:covering}.
We further use the existence of geometric models for the Bianchi groups in which all edge stabilisers are finite cyclic
and all cells of dimension $2$ and higher are trivially stabilised.
Therefore, the system of finite subgroups of the Bianchi groups admits inclusions only emanating from cyclic groups.
This makes the Bianchi groups and their subgroups subject to the splitting of Bredon homology stated in Theorem~\ref{splitting}.

The proof of Theorem~\ref{splitting} is based on the above particularities of the Bianchi groups,
and applies the following splitting lemma for the involved representation rings
 to a Bredon complex for~$(\underline{\rm E}\Gamma, \Gamma)$.

\vbox{\begin{lemma}[\cite{Rahm:equivariant}] \label{representation ring splitting}
Consider a group $\Gamma$ such that every one of its finite subgroups is either cyclic of order at most~$3$, or of one of the types
$\Dtwo, \Dthree$ or~$\Afour$.
Then there exist bases of the complex representation rings of the finite subgroups of~$\Gamma$,
 such that simultaneously every morphism of representation rings
 induced by inclusion of cyclic groups into finite subgroups of~$\Gamma$,
 splits as a matrix into the following diagonal blocks.
\begin{enumerate}
 \item A block of rank $1$ induced by the trivial and regular representations,
 \item a block induced by the $2$--torsion subgroups
 \item and a block induced by the $3$--torsion subgroups.
\end{enumerate}
\end{lemma}
}

As this splitting holds simultaneously for every morphism of representation rings,
we have such a splitting for every morphism of formal sums of representation rings,
and hence for the differential maps of the Bredon complex for any Bianchi group and any of their subgroups.

The bases that are mentioned in the above lemma, are obtained by elementary base transformations
from the canonical basis of the complex representation ring of a finite group to
a basis whose matrix form has
\begin{itemize}
 \item its first row concentrated in its first entry, for a finite cyclic group (edge stabiliser).
The base transformation is carried out by summing over all representations to replace the trivial representation by the regular representation.
 \item its first column concentrated in its first entry, for a finite non-cyclic group (vertex stabiliser).
The base transformation is carried out by subtracting the trivial representation from each representation, except from itself.
 \end{itemize}
The details are provided in \cite{Rahm:equivariant}.

In this setting, the technique has inspired work beyond the range of arithmetic groups,
which has led to formulas for the integral Bredon homology and equivariant K-homology of all compact 3-dimensional hyperbolic reflection groups~\cite{LORS},
through a novel criterion for torsion-freeness of equivariant K-homology in a more general framework.

\subsection{Chen--Ruan orbifold cohomology of the complexified Bianchi orbifolds} \label{orbifold state}

The action of the Bianchi groups $\Gamma$ on real hyperbolic $3$-space $\symmetricspace = \mathrm{SL}_2(\C)/\mathrm{SU}_2$
induces an action on a complexification $\symmetricspace_\C$ of the latter (of real dimension $6$).
For the orbifolds $[\symmetricspace_\C /_\Gamma]$ given by this action, we can compute the Chen--Ruan Orbifold Cohomology as follows.

Let~$\Gamma$ be a discrete group acting \emph{properly}, i.e. with finite stabilizers, by diffeomorphisms on a manifold~$X$.
For any element $g \in \Gamma$, denote by $C_\Gamma(g)$ the centralizer of $g$ in~$\Gamma$. Denote by $X^g$ the subset of $X$ consisting of the fixed points of $g$.

\begin{df}
Let $T \subset \Gamma$ be a set of representatives of the conjugacy classes of elements of finite order in~$\Gamma$.
Then we set
$$ \Homol^*_{orb}([X / _\Gamma]) := \bigoplus_{g \in T} \Homol^* \left( X^g / C_\Gamma(g); \thinspace \rationals \right),$$
where $\Homol^* \left( X^g / C_\Gamma(g); \thinspace \rationals \right)$ is the ordinary cohomology of the quotient space $X^g / C_\Gamma(g)$.
\end{df}
It can be checked that this definition gives the vector space structure of the orbifold cohomology defined by Chen and Ruan~\cite{ChenRuan}, if we forget the grading of the latter.
We can verify this fact using arguments analogous to those used by Fantechi and G\"ottsche \cite{FantechiGoettsche}
in the case of a finite group~$\Gamma$ acting on~$X$.
The additional argument needed when considering some element $g$ in~$\Gamma$ of infinite order, is the following.
As the action of~$\Gamma$ on $X$ is proper, $g$ does not admit any fixed point in $X$. Thus,
$ \Homol^* \left( X^g / C_\Gamma(g); \thinspace \rationals \right) = \Homol^* \left( \emptyset; \thinspace \rationals \right) = 0. $

\begin{theorem}[\cite{PerroniRahm}] \label{introduced result}
Let $\Gamma$ be a finite index subgroup in a Bianchi group (except over the Gaussian or Eisensteinian integers).
Denote by $\lambda_{2n}$ the number of conjugacy classes of cyclic subgroups of order ${n}$
in $\Gamma$.
Denote by $\lambda_{2n}^*$
the cardinality of the subset of conjugacy classes which are contained in a dihedral subgroup of order $2n$
in~$\Gamma$. Then,
\small
$$ \Homol^d_{orb}\left([\symmetricspace_\C /_\Gamma] \right) \cong
\Homol^d\left(\symmetricspace/_{\Gamma}; \thinspace \rationals \right) \oplus
\begin{cases}
\rationals^{\lambda_4 +2\lambda_6 -\lambda_6^*}, & d=2, \\
 \rationals^{\lambda_4-\lambda_4^* +2\lambda_6 -\lambda_6^*}, & d=3, \\
0, & \mathrm{otherwise}. \end{cases}$$
\normalsize
\end{theorem}

The (co)homology of the quotient space $\symmetricspace  / _\Gamma$
 has been computed numerically for a large range of Bianchi groups \cite{Vogtmann}, \cite{Scheutzow}, \cite{Rahm:higher_torsion};
 and bounds for its Betti numbers are given in \cite{Kraemer:Thesis}.
Kr\"amer \cite{Kraemer:Diplom} has determined number-theoretic formulas
 for the numbers $\lambda_{2n}$ and $\lambda_{2n}^*$ of conjugacy classes of finite subgroups in the Bianchi groups.

Building on this, Perroni and the author established the following result~\cite{PerroniRahm}.
\begin{theorem}\label{mainthm} ${}$\\
Let $\symmetricspace_\C /_\Gamma$ be the coarse moduli space of
$[\symmetricspace_\C /_\Gamma]$;  \\
and let $Y \ra \symmetricspace_\C /_\Gamma$ be a crepant resolution of
$\symmetricspace_\C /_\Gamma$.\\
Then there is an isomorphism as graded $\CC$-algebras between the Chen-Ruan
cohomology ring of $[\symmetricspace_\C /_\Gamma]$ and the singular cohomology ring of $Y$:
$$
\left( H_{\rm CR}^*([\symmetricspace_\C / _\Gamma]) , \cup_{\rm CR} \right) \cong
\left( H^*(Y) , \cup \right) \, .
$$
\end{theorem}
The Chen--Ruan orbifold cohomology is conjectured by Ruan to match the quantum corrected classical cohomology ring of a crepant resolution for the orbifold.
Perroni and the author proved furthermore that the Gromov-Witten invariants involved in the definition
of the quantum corrected cohomology ring of $Y\ra \symmetricspace_\C /_\Gamma$
vanish. Hence, they deduced the following.
\begin{corollary}\label{maincor}
Ruan's crepant resolution conjecture holds true for the complexified Bianchi orbifolds $[\symmetricspace_\C / _\Gamma]$.
\end{corollary}

A result on the vector space structure of the Chen--Ruan orbifold cohomology of Bianchi orbifolds are the below two theorems.

\begin{theorem}[\cite{Rahm:equivariant}] \label{3-torsion quotients}
For any element $\gamma$ of order $3$ in a finite index subgroup~$\Gamma$ in a Bianchi group with units~$\{\pm 1\}$,
 the quotient space $\Hy^\gamma /_{C_\Gamma(\gamma)}$ of the rotation axis modulo the centralizer of $\gamma$
 is homeomorphic to a circle.
\end{theorem}

\begin{theorem}[\cite{Rahm:equivariant}] \label{2-torsion quotients}
Let $\gamma$ be an element of order $2$ in a Bianchi group~$\Gamma$ with units~$\{\pm 1\}$.
Then, the homeomorphism type of the quotient space $\Hy^\gamma /_{C_\Gamma(\gamma)}$ is
\end{theorem}
\begin{itemize}
\item[$\edgegraph$]
\textit{an edge without identifications, if $\langle \gamma \rangle$ is contained in a subgroup of type $\Dtwo$
inside~$\Gamma$ and
\item[$\circlegraph$]
a circle, otherwise.}
\end{itemize}

Denote by $\lambda_{2\ell}$ the number of conjugacy classes of subgroups of type $\Z/\ell$
in a finite index subgroup {$\Gamma$} in a Bianchi group with units $\{\pm 1 \}$.
Denote by $\lambda_{2\ell}^*$ the number of conjugacy classes of subgroups of type $\Z/\ell$
which are contained in a subgroup of type $\Dtwon$ in~$\Gamma$.
By \cite{Rahm:equivariant},
there are \mbox{$2\lambda_6 -\lambda_6^*$} conjugacy classes of elements of order $3$.
As a result of Theorems~\ref{3-torsion quotients} and~\ref{2-torsion quotients},
 the vector space structure of the orbifold cohomology of $[\Hy_\R / _\Gamma]$ is given as

 $
\Homol^\bullet_{orb}([\Hy_\R / _\Gamma]) \cong
\Homol^{\bullet} \left( \Hy_\R / _\Gamma; \thinspace \rationals \right)
\bigoplus\nolimits^{\lambda_4^*} \Homol^{\bullet} \left( \edgegraph; \thinspace \rationals \right)
\bigoplus\nolimits^{(\lambda_4 -\lambda_4^*)} \Homol^{\bullet} \left( \circlegraph; \thinspace \rationals \right)
\bigoplus\nolimits^{(2\lambda_6 -\lambda_6^*)} \Homol^{\bullet} \left(\circlegraph; \thinspace \rationals \right).
$
\normalsize \\
The (co)homology of the quotient space $\Hy_\R  / _\Gamma$
 has been computed numerically for a large range of Bianchi groups \cite{Vogtmann}, \cite{Scheutzow}, \cite{Rahm:higher_torsion};
 and bounds for its Betti numbers are given in \cite{Kraemer:Thesis}.
Kr\"amer \cite{Kraemer:Diplom} has determined number-theoretic formulas
 for the numbers $\lambda_{2\ell}$ and $\lambda_{2\ell}^*$ of conjugacy classes of finite subgroups in the full Bianchi groups.
Kr\"amer's formulas were evaluated for hundreds of thousands of Bianchi groups \cite{Rahm:formulas},
and these values are matching with the ones from the orbifold structure computations with \cite{Rahm:BianchiGP}
in the cases where the latter are available.

When we pass to the complexified orbifold $[\Hy_\C / _\Gamma]$,
the real line that is the rotation axis in~$\Hy_\R$ of an element of finite order, becomes a complex line.
However, the centralizer still acts in the same way by reflections and translations.
So, the interval $\edgegraph$ as a quotient of the real line yields a stripe
 $\edgegraph \times \R$ as a quotient of the complex line.
And the circle $\circlegraph$ as a quotient of the real line yields a cylinder
 $\circlegraph \times \R$ as a quotient of the complex line.
Therefore,  using the degree shifting numbers computed in~\cite{Rahm:equivariant},
we obtain the result of Theorem~\ref{introduced result}.

\subsection*{Acknowledgement.} The author was supported by the MELODIA project, grant number ANR-20-CE40-0013, during the revision of this paper.

\end{document}